\documentclass[12pt]{article}
\usepackage{amsmath, amssymb, latexsym}
\renewcommand{\baselinestretch}{1.15}
 \newtheorem{theorem}{Theorem}[section]
\newtheorem{lemma}[theorem]{Lemma}
\newtheorem{proposition}[theorem]{Proposition}

\newtheorem{remark}[theorem]{Remark}

\newtheorem{definition}[theorem]{Definition}

\renewcommand{\title}[1]
{\thispagestyle{empty}
\begin{center}
{\Large \bf #1}
\end{center}}

\newcommand{\authors}[1]
{\begin{center}
\renewcommand{\thefootnote}{\fnsymbol{footnote}}
\setcounter{footnote}{3} {\sc #1 }
\end{center}}

\newcommand{\ack}[1]{\footnote{#1}}

\newcommand{\address}[1]
{\vskip 5ex
\renewcommand{\baselinestretch}{1}
\footnotesize \normalsize
 #1 \\
}

\begin{document}

\title{Mixed norm and multidimensional Lorentz spaces}

\authors{
Sorina Barza\ack{Research partially supported by KAW 2000.0048 and STINT KU
2002-4025.},
Anna Kami\'nska,  
Lars-Erik Persson,  
and
Javier Soria\ack{Research partially supported by Grants BFM2001-3395,   
2001SGR00069 and The Swedish Research Council no. 624-2003-571.} }

\bigskip

 {\narrower\noindent \textbf{Abstract.} \small{In the last decade, the problem of characterizing the normability of the weighted Lorentz spaces has been completely solved (\cite{Sa}, \cite{CaSoA}). However, the question for multidimensional Lorentz spaces is still open. In this paper, we consider weights of product type, and  give necessary and sufficient conditions for
the Lorentz spaces, defined with respect to the two-dimensional decreasing
rearrangement, to be normable. To this end, it is also useful to study the mixed norm Lorentz spaces. Finally, we prove embeddings between all the classical, multidimensional, and mixed norm Lorentz spaces.}\par}

\bigskip

\section{Introduction}

Let $f$ $:\mathbb{R}^{n}\rightarrow \mathbb{R}$ be a Lebesgue measurable
function. The usual decreasing rearrangement of $f$ on $(0,\infty )$ is
given by $\ \ f^{\ast }(t)=\inf \{\sigma :\lambda _{f}(\sigma )\leq t\},$ $
t>0,$ where $\lambda _{f}(\sigma )=\left| \{x:\left| f(x)\right| >\sigma
\}\right| $ is the distribution function (see e.g. \cite{BS}). In \cite
{BaPeSo}, a multidimensional decreasing rearrangement  was defined by using
the \lq\lq Layer cake formula,'' which recovers a function by means of its level
sets. Surprisingly,   this definition coincides with the
multivariate rearrangement defined in \cite{Bl}, which we will use in this paper. For simplicity we are going
to reduce our definitions to the two-dimensional case because the extensions
to higher dimensions only require natural modifications. By $f_{y}^{\ast }(x,t)$ we will denote the
decreasing rearrangement of $f$ with respect to the second variable $y$,
under fixed first variable $x$, and $f_{x}^{\ast }(s,y)$ will denote the
decreasing rearrangement of $f$ with respect to the first variable $x$, under
fixed second variable $y.$ The multivariate decreasing rearrangement of $f$,
  first with respect to the second variable $y$ and then with respect to $x$, will be denoted by $ 
f_{2,1}^{\ast }(s,t)=f_{yx}^{\ast }(s,t)=(f_{y}^{\ast }(\cdot ,t))_{x}^{\ast
}(s)$, and similarly,  $f_{1,2}^{\ast }(s,t)=f_{xy}^{\ast }(s,t):=(f_{x}^{\ast }(s,\cdot
))_{y}^{\ast }(t).$ Associated with a function $f(x,y)$ we define the
following   multivariate averaging operators of Hardy type (see \cite{KP} for more information):

\begin{eqnarray*}
S^{2}f(s,t)&:=&\frac{1}{st}\int_{0}^{s}\int_{0}^{t}f(\sigma,\tau)\,d\tau\,d\sigma,
 \\
f^{\ast \ast }(s,t)&:=&S^{2}f_{yx}^{\ast }(s,t)=\frac{1}{st}
\int_{0}^{s}\int_{0}^{t}f_{yx}^{\ast }(\sigma ,\tau )\,d\sigma \,d\tau,
\\
S_{2,1}f(s,t)&:=&\frac{1}{s}\int_{0}^{s}\left( \frac{1}{ 
t}\int_{0}^{t}f_{y}^{\ast }(\cdot ,\tau )\,d\tau \right) _{x}^{\ast }(\sigma
)\,d\sigma.
\end{eqnarray*}
We will also  use the notations $f_{x}^{\ast \ast }(s,y)=\frac{1}{s} 
\int_{0}^{s}f_{x}^{\ast }(\sigma ,y)\,d\sigma $ and $f_{y}^{\ast \ast }(x,t)= 
\frac{1}{t}\int_{0}^{t}f_{y}^{\ast }(x,v)dv$ and hence,
$S_{2,1}f(s,t)=(f_{y}^{\ast \ast }(\cdot ,t))_{x}^{\ast \ast }( s)$ (for this reason, we will sometimes write $f_{yx}^{\ast \ast }(s,t):=S_{2,1}f(s,t)$.)

We recall the definition of the classical Lorentz space: If $v$ is a weight in $ 
\mathbb{R}_{+}$, that is, $v$ is non-negative and locally integrable,  and $0<p<\infty ,$ then
\begin{equation}
\Lambda ^{p}(\mathbb{R}^n,v)=\left\{ f:\mathbb{R}^{n}\rightarrow \mathbb{R};\ \left\|
f\right\| _{\Lambda ^{p}(v)}:=\bigg(\int_{0}^{\infty }(f^{\ast
}(t))^{p}v(t)\,dt\bigg)^{1/p}<\infty \right\} .  \label{Lor}
\end{equation}
When $n=1$, we shall write $\Lambda^p(v)=\Lambda^p(\mathbb{R},v)$.

Below we give two different definitions of two-dimensional Lorentz spaces. The first one is based on a concept of \lq\lq mixed norms", while the second one is based on classical definitions and multivariate decreasing rearrangement $f^*_{yx}$ or $f^*_{xy}$. We shall see later on that these spaces are essentially different.

According to \cite{Bl}, if $u$ and $v$ are weights in $\mathbb{R}_{+}$ and $ 
0<p,q<\infty ,$ we say that a measurable function $f$ belongs to the mixed weighted Lorentz space $\Lambda^{q} ({u})[\Lambda^{p} ({v})]$, if 
\begin{equation*}
 \left\| f\right\| _{\Lambda^{q} ({u})[\Lambda^{p} ({v})]}:= 
\left( \int_{0}^{\infty }\bigg[\left( \int_{0}^{\infty }(f_{y}^{\ast}(\cdot
,t))^pv(t)\,dt\right) _{x}^{\ast }(s)\bigg]^{q/p}u(s)\,ds\right) ^{1/q}<\infty.
\end{equation*}
Similarly, we  say that   $f$   belongs to the
two-dimensional Lorentz space $\Lambda_{2}^{p}(w),$ provided 

\begin{equation}
\bigskip\left\| f\right\|_{\Lambda_{2}^{p}(w)}
:=\bigg(\int_{\mathbb R^2_+}(f_{yx}^{\ast}(s,t))^{p}w(s,t)\,ds\,dt\bigg)^{1/p}<\infty,
\label{Lor2dim}
\end{equation}
 where $w$ is a weight
function defined on $\mathbb{R}^{2}_+$ (see \cite{BaPeSo}).

In Section \ref{nols}   we prove the most important properties of the rearrangements $
f_{yx}^{\ast }(s,t),$ $f^{\ast \ast }(s,t),$ and  $f_{yx}^{\ast \ast }(s,t)$. We show, among other things, that $f^{**}_{yx}(s,t)$ is a sublinear operator, analogously as the Hardy operator $t^{-1}\int_0^t f^*(s)\,ds$ in the case of one variable, while $f^{**}(s,t)$ does not enjoy this property. We also
give equivalent conditions for the normability of $\Lambda _{2}^{p}(w)$, if $ 
w(s,t)=u(s)v(t)$, as well as the space $\Lambda^{p} ({u})[\Lambda^{p} ({v})]$. The first results in this theory are due to Lorentz (see \cite{Lo}), where he characterized when   the functional defined in (\ref{Lor}) is a norm on $\Lambda^p(v)$  (similar results for the spaces $\Lambda_{2}^{p}(w)$ have been  recently proved in \cite{BaPeSo}). In  \cite{Sa}, Sawyer extended Lorentz result to characterize  the normability of $\Lambda^p(v)$, for $p>1$, and the case $p=1$ was established in \cite{CaSoA}.  In Section \ref{nols} we also compare the spaces $\Lambda_2^p(uv)$ and $\Lambda^{p} ({u})[\Lambda^{p} ({v})]$, showing in particular that they do not coincide. In Section \ref{emth}, we show the embeddings between some of  the
spaces defined above. Our theorems generalize
previous results from \cite{Bl} and \cite{Y}, where they considered the case of power weights; i.e., when the Lorentz space is of the form $L^{p,q}$, $0<p,q<\infty$ (see \cite{BS}).

\section{ Normability of two-dimensional Lorentz spaces}\label{nols}
The main technique to prove normability of the classical   Lorentz spaces is given in terms of the boundedness of the Hardy operator for the class of monotone functions, and the fact that this transformation enjoys a subadditive property. We will show that in higher dimensions, this operator has to be replaced by $f^{**}_{yx}$, since the natural generalization $f^{**}(s,t)$ of the Hardy operator is not sublinear. First, we will describe some properties of the transformations  $f^{\ast\ast}(s,t)$
and $f_{yx}^{\ast\ast}(s,t),$ which will be used later on.

\begin{proposition}
\label{properties}
With the notations above we have: 
\begin{enumerate}
\item
$f_{yx}^{\ast }(s,t)\leq f_{yx}^{\ast \ast }(s,t).$

\item
$(f+g)_{yx}^{\ast \ast }(s,t)\leq f_{yx}^{\ast \ast }(s,t)+g_{yx}^{\ast \ast
}(s,t).$
\item
$f_{yx}^{\ast \ast }(s,t)\leq f^{\ast \ast }(s,t),$
 (in general $ $ $f_{yx}^{\ast \ast }(s,t)\not=f^{\ast \ast }(s,t)$.)

\item  
 $f^{\ast \ast }$ is not sublinear, but we have:
\begin{equation*}
(f+g)^{\ast \ast }(s,t)\leq 4(f^{\ast \ast }(s,t)+g^{\ast \ast }(s,t)).
\end{equation*}
\end{enumerate}
\end{proposition}

\noindent{\bf Proof:} 
(1.) By classical arguments we have $\ \frac{1}{t}\int_{0}^{t}f_{y}^{\ast
}(x,\tau )\,d\tau \geq f_{y}^{\ast }(x,t)$, for any $t$ and $x.$
Hence 
\begin{equation*}
\bigg(\frac{1}{t}\int_{0}^{t}f_{y}^{\ast }(\cdot ,\tau )\,d\tau \bigg)_{x}^{\ast
}(\sigma )\geq (f_{y}^{\ast }(\cdot ,t))_{x}^{\ast }(\sigma )=f_{yx}^{\ast
}(\sigma,t ).
\end{equation*}
Integrating now with respect to $\sigma $ we get

\begin{equation*}
f_{yx}^{\ast \ast }(s,t)=\frac{1}{s}\int_{0}^{s}\left( \frac{1}{t} 
\int_{0}^{t}f_{y}^{\ast }(\cdot ,\tau )\,d\tau \right) _{x}^{\ast }(\sigma
)\,d\sigma \geq \frac{1}{s}\int_{0}^{s}f_{yx}^{\ast }(\sigma ,t )\,d\sigma \geq
f_{yx}^{\ast }(s,t),
\end{equation*}
and the inequality is proved.

(2.) We use now the following fact (see  \cite{KPS}): if $v$ is a decreasing function, then 
\begin{equation}
\sup_{\rho }\int_{\mathbb{R}}|f(x)|v(\rho (x))\,dx=\int_{0}^{\infty }f^{\ast
}(t)v(t)\,dt,  \label{f1}
\end{equation}
where the supremum is taken over all measure preserving transformations $\rho
:\mathbb{R}\rightarrow \mathbb{R}_{+}.$ Therefore, 
using (\ref{f1}) and the subadditivity of the one-dimensional maximal function:
\begin{align*}
(f+g)_{yx}^{\ast \ast }(s,t)& =\frac{1}{s}\int_{0}^{s}\left( \frac{1}{t} 
\int_{0}^{t}(f+g)_{y}^{\ast }(\cdot , \tau)\,d\tau \right) _{x}^{\ast }(\sigma
)\,d\sigma  \\
& =\frac{1}{s}\int_{0}^{s}\left( \sup_{\rho }\int_{\mathbb{R}}|f+g|(\cdot ,y) 
\frac{1}{t}\chi _{(0,t)}(\rho (y))\,dy\right) _{x}^{\ast }(\sigma )\,d\sigma \\
& \leq \frac{1}{s}\int_{0}^{s}\left( \sup_{\rho }\int_{\mathbb{R}}|f(\cdot ,y) |
\frac{1}{t}\chi _{(0,t)}(\rho (y))\,dy\right.   \\
& \qquad+\left. \sup_{\rho }\int_{\mathbb{R}}|g(\cdot ,y)|\frac{1}{t}\chi
_{(0,t)}(\rho (y))\,dy\right) _{x}^{\ast }(\sigma )\,d\sigma  \\
& \leq \frac{1}{s}\int_{0}^{s}\left( \sup_{\rho }\int_{\mathbb{R}}^{\
}|f(\cdot ,y)|\frac{1}{t}\chi _{(0,t)}(\rho (y))\,dy\right) _{x}^{\ast }(\sigma
)\,d\sigma  \\
&  \qquad+\frac{1}{s}\int_{0}^{s}\left( \sup_{\rho }\int_{\mathbb{R}}|g(\cdot ,y) |
\frac{1}{t}\chi _{(0,t)}(\rho (y))\,dy\right) _{x}^{\ast }(\sigma )\,d\sigma  \\
& =\frac{1}{s}\int_{0}^{s}\left( \frac{1}{t}\int_{0}^{t}f_{y}^{\ast }(\cdot
,\tau )\,d\tau \right) _{x}^{\ast }(\sigma )\,d\sigma  \\
& \qquad +\frac{1}{s}\int_{0}^{s}\left( \frac{1}{t}\int_{0}^{t}g_{y}^{\ast }(\cdot
,\tau )\,d\tau \right) _{x}^{\ast }(\sigma )\,d\sigma  \\
& =f_{yx}^{\ast \ast }(s,t)+g_{yx}^{\ast \ast }(s,t),
\end{align*}
which completes the proof.

(3.) By\ using (\ref{f1}) we get,

\begin{align*}
f^{**}_{yx}(s,t)&=  
\frac{1}{s}\int_{0}^{s}\left( \frac{1}{t}\int_{0}^{t}f_{y}^{\ast }(\cdot
,\tau )\,d\tau \right) _{x}^{\ast }(\sigma )\,d\sigma \\ & =\frac{1}{st}\sup_{\rho
}\int_{\mathbb{R}}\left( \int_{0}^{t}f_{y}^{\ast }(x,\tau )\,d\tau \right)
\chi _{\lbrack 0,s]}(\rho (x))\,dx \\
& =\frac{1}{st}\sup_{\rho }\int_{0}^{t}\bigg(\int_{\mathbb{R}}f_{y}^{\ast
}(x,\tau )\chi _{\lbrack 0,s]}(\rho (x))\,dx\bigg)\,d\tau  \\
& \leq \frac{1}{st}\int_{0}^{t}\left( \sup_{\rho }\int_{\mathbb{R} 
}f_{y}^{\ast }(x,\tau )\chi _{\lbrack 0,s]}(\rho (x))\,dx\right) \,d\tau  \\
& =\frac{1}{st}\int_{0}^{t}\left( \int_{0}^{s}(f_{y}^{\ast }(\cdot ,\tau
))_{x}^{\ast }(\sigma)\,d\sigma \right) \,d\tau  \\
& =\frac{1}{st}\int_{0}^{s}\int_{0}^{t}f_{yx}^{\ast }(\sigma ,\tau )\,d\tau
\,d\sigma  \\
& =f^{\ast \ast }(s,t).
\end{align*}
In general,  $f^{\ast \ast }(s,t)\neq $ $f_{yx}^{\ast \ast }(s,t).$ To see this, it suffices to consider the function  $f=\chi _{D},$ where $D=[0,3]\times \lbrack 0,1]\cup \lbrack 2,3]\times
\lbrack 1,2].$ 

(4.)  The first part of the statement is a consequence of \cite[Theorem 3.7]
{BaPeSo} (just observe that the weight $(st)^{-1}\chi_{[0,s]}(\sigma)\chi_{[0,t]}(\tau)$ does not satisfy the hypothesis of this theorem). Using now, e.g. \cite[Lemma 2.2 IV]{Bl}  we find that
\begin{align*}
\left( f+g\right) ^{\ast \ast }(s,t)& =\frac{1}{st}\int_{0}^{s} 
\int_{0}^{t}(f+g)_{yx}^{\ast }(\sigma,\tau)\,d\tau\,d\sigma \\
& =\frac{1}{st}\int_{0}^{\infty }\int_{0}^{\infty }(f+g)_{yx}^{\ast
}(\sigma,\tau)\chi _{\lbrack 0,s]}(\sigma)\chi _{\lbrack 0,t]}(\tau)\,d\tau\,d\sigma\\
& \leq \frac{1}{st}\int_{0}^{\infty }\int_{0}^{\infty }f_{yx}^{\ast }\Big(\frac{\sigma 
}{2},\frac{\tau}{2}\Big)\chi _{\lbrack 0,s]}(\sigma)\chi _{\lbrack 0,t]}(\tau)\,d\tau\,d\sigma\\
& \qquad+\frac{1}{st}\int_{0}^{\infty }\int_{0}^{\infty }g_{yx}^{\ast }\Big(\frac{\sigma}{2} 
,\frac{\tau}{2}\Big)\chi _{\lbrack 0,s]}(\sigma)\chi _{\lbrack 0,t]}(\tau)\,d\tau\,d\sigma\\
& =\frac{4}{st}\int_{0}^{\infty }\int_{0}^{\infty }f_{yx}^{\ast }(\sigma,\tau)\chi
_{\lbrack 0,s]}(2\sigma)\chi _{\lbrack 0,t]}(2\tau)\,d\tau\,d\sigma\\
&\qquad+ \frac{4}{st}\int_{0}^{\infty }\int_{0}^{\infty }g_{yx}^{\ast }(\sigma,\tau)\chi
_{\lbrack 0,s]}(2\sigma)\chi _{\lbrack 0,t]}(2\tau)\,d\tau\,d\sigma\\
& \leq \frac{4}{st}\int_{0}^{\infty }\int_{0}^{\infty }f_{yx}^{\ast
}(\sigma,\tau)\chi _{\lbrack 0,s]}(\sigma)\chi _{\lbrack 0,t]}(\tau)\,d\tau\,d\sigma\\
& \qquad+\frac{4}{st}\int_{0}^{\infty }\int_{0}^{\infty }g_{yx}^{\ast }(\sigma,\tau)\chi
_{\lbrack 0,s]}(\sigma)\chi _{\lbrack 0,t]}(\tau)\,d\tau\,d\sigma\\
& =\frac{4}{st}\int_{0}^{s}\int_{0}^{t}f_{yx}^{\ast }(\sigma,\tau)\,d\tau\,d\sigma+\frac{4}{st} 
\int_{0}^{s}\int_{0}^{t}g_{yx}^{\ast }(\sigma,\tau)\,d\tau\,d\sigma,
\end{align*}
which completes the proof of the last statement.
$\hfill\Box$\bigbreak
$\ $

The problem of finding conditions on $v$ such that $\Lambda ^{p}(v), $ defined
in (\ref{Lor}), is normable (in fact, a Banach space, since completeness always holds), was solved for $ 
p>1$, by E. Sawyer (\cite{Sa}). This condition is that the Hardy-Littlewood
maximal operator is bounded on $\Lambda ^{p}(v).$ The weights for which this
holds were first characterized by M. A. Ari\~{n}o and B. Muckenhoupt \cite
{AM}, and it is known as the $B_{p}$ condition: there exists $C>0$ such that,
for all $r>0,$ 
\begin{equation}\label{bcondition}
r^{p}\int_{r}^{\infty }\frac{v(x)}{x^{p}}\,dx\leq C\int_{0}^{r}v(x)\,dx.
\end{equation}

It is clear that (\ref{bcondition}) i\vspace{0in}s not the right condition
for $p=1$, since with $v\equiv 1 $ we have that\ $\Lambda ^{1}(v)=L^{1}$, which
is a Banach space, but $v$ does not satisfy (\ref{bcondition}). This endpoint case
was solved by M.\  J.\  Carro, A. Garc\'{\i}a del Amo and J.\ Soria in \cite{CaSoA}. Now,  the weight
has to satisfy the so called $B_{1,\infty}$ condition: there exists $C>0$ such
that for all $0<s\leq r<\infty ,$ 
\begin{equation}
\frac{1}{r}\int_{0}^{r}v(x)\,dx\leq \frac{C}{s}\int_{0}^{s}v(x)\,dx.
\label{bpcondition}
\end{equation}
This motivates the consideration of the same type of problems for the
two-dimensional Lorentz space $\Lambda _{2}^{p}(w).$ The characterization of the weights $w$ such that $\left\| f\right\|_{\Lambda _{2}^{p}(w)} $ is a norm  was proved in  \cite{BaPeSo} (the corresponding result for $\Lambda^p(v)$ was proved in \cite{Lo}), and there it was
also shown that if $\Lambda _{2}^{p}(w)$ is a Banach space, then $p\geq 1.$

The normability conditions for the space $\Lambda^{p}(u)[\Lambda^{p}(v)]$
follow from the general theory of mixed norm spaces. Let $(\Omega ,\Sigma ,\mu )$ be a $\sigma $-finite measure space, where $\Sigma$ is a $\sigma$-algebra of subsets of $\Omega$ and $\mu$ is a $\sigma$-finite measure on $\Sigma$. Letting $L^0(\mu)$ be the space of all real-valued $\Sigma$-measurable functions on $\Omega$, the space $E\subset L^0(\Sigma)$ equipped with a quasi-norm $||\cdot||_E$ is called a quasi-normed function lattice if for $f\in L^0(\mu)$, $g\in E$, and $|f|\le|g|$ a.e., we have that $f\in E$, and $||f||_E\le||g||_E$. Given two quasi-Banach function lattices $E$ and $F$ defined on $(\Omega _{1},\Sigma _{1},\mu _{1})$ and $(\Omega_{2},\Sigma _{2},\mu _{2})$, respectively, the mixed quasi-normed space $E[F]$ consists of all $\Sigma _{1}\times\Sigma _{2}$-measurable functions $f:\Omega _{1}\times\Omega _{2}\rightarrow \mathbb{R}$ such that $||f||_{E[F]}:=||x\rightarrow||f(x,\cdot)||_F||_E<\infty.$ We have the following result:

\begin{proposition}
\label{Lmix}
A mixed quasi-normed space $E[F]$ is normable, if and only if both $E$ and $F$ are normable. In particular,  if $p>1,$ then $\Lambda^{p}(u)[\Lambda^{p}(v)]$ is a Banach space if
and only if $u$ and $v\in B_{p}$, and if \ $p=1,$ $\Lambda^{1}(u)[\Lambda
^{1}(v)]$ is a Banach space if and only if $u$ and $v\in B_{1,\infty}$.
\end{proposition}

\noindent{\bf Proof:} 
Assume first that $E$ and $F$ are normable, i.e. there exist norms $|||\cdot|||_E$ on $ 
E,$   $|||\cdot|||_F$ on $F$, such that $A_{1}\left\| \cdot \right\| _{E}\leq ||| \cdot|||_{E} \leq B_{1}\left\| \cdot \right\| _{E}$
and $A_{2}\left\| \cdot \right\| _{F}\leq ||| \cdot ||| _{F} \leq B_{2}\left\| \cdot \right\| _{F}$ \ for some $A_j,B_j>0,\ j=1,2.$  We define $||| f||| :=||| x\rightarrow |||
f(x,\cdot )||| _{F}||| _{E}$ and we want to show
that $||| \cdot |||$ is a norm which is
equivalent to $\left\| f\right\| _{E[F]}.$ We have 
\begin{align*}
||| f+g|||& =|||x\rightarrow |||f(x,\cdot )+g(x,\cdot )||| _{F}|||_{E} \\
& \leq ||| x\rightarrow||| f(x,\cdot )|||_{F}+x\rightarrow||| g(x,\cdot )|||
_{F}|||_{E} \\
& \leq ||| x\rightarrow ||| f(x,\cdot )||| _{F}||| _{E}+|||x\rightarrow |||
g(x,\cdot )|||_{F}||| _{E} \\
& =|||f|||+||| g|||.
\end{align*}
On the other hand, since the norm on $E$ has the lattice property, it is clear that 
$$
A_1A_2||x\rightarrow||f(x,\cdot)||_F||_E\le||| f||| \le B_1B_2||x\rightarrow||f(x,\cdot)||_F||_E,
$$
which shows that $E[F]$ is normable. 

Conversely, suppose now that $E[F]$ is a Banach function space, i.e. there
exists a lattice norm $||| \cdot ||| $ which is
equivalent to $\left\| \cdot \right\| _{E[F]}.$ If $Tf(s,t)=f(s)g(t),$
where $f$ is an arbitrary function in $E$ and $g$ is a fixed function in $F$
such that $\left\| g\right\| _{F}=1,$ then we have that $Tf\in E[F]$ and $\left\|
Tf\right\| _{E[F]}=\left\| f\right\| _{E},$ which shows that $T$ is an isometric embedding of $E$ into $E[F]$. Define $|||f||| _{E}:=|||Tf||| .$ It is obvious that $||| f||| _{E}$ is equivalent to $\|f\|_{E}, $ and  also
$$
|||f_{1}+f_{2}||| _{E}\le ||| f_{1}||| _{E}+||| f_{2}||| _{E.}
$$
Hence $|||\cdot |||_{E}$ is a norm and $\ E$ is
normable. Similarly one can prove that $F$ is normable. The second part of the statement is a consequence of  the results in \cite{Sa} and \cite{CaSoA}, respectively.
$\hfill\Box$\bigbreak

We will
now prove the main theorem of this Section. We characterize the normability
of\ the two-dimensional Lorentz space $\Lambda_{2}^{p}(w)$ \vspace{0in}in
the particular case when $w(s,t)=u(s)v(t).$ In this case we will use the
notation $\Lambda_{2}^{p}(uv).$ By $\ L_{\rm{dec}}^{p}(\mathbb{R}^{2}_{+},uv)$ we
denote the cone of functions in $L^{p}(\mathbb{R}^{2}_{+},uv)$, decreasing in each variable separately. We will need first the following simple lemma:

\begin{lemma}
\label{LH}The two-dimensional Hardy operator $S^{2}$ is bounded from $ 
L_{\rm{dec}}^{p}(\mathbb{R}^{2}_{+},w)$ to $\Lambda_2^p(w)$ if
and only if \ $S_{2,1}$ is bounded from  $\Lambda _{2}^{p}(w)$ to $\Lambda _{2}^{p}(w)$ .
\end{lemma}

\noindent{\bf Proof:} 
It is easy to see that on the cone of decreasing functions in each variable
separately the operators $S^{2}$ and $S_{2,1}$ coincide.$\hfill\Box$\bigbreak

\begin{theorem}\label{mainth}
Let $p>1.$ The following conditions are equivalent:

\begin{enumerate}
\item  $\Lambda _{2}^{p}(uv)$ \textit{is normable.}

\smallskip

\item  $u$ \textit{and} $\ v$ $\in B_{p}.$

\smallskip

\item  $\Lambda^{p} ({u})[\Lambda^{p} ({v})]$ \textit{is normable.} 
\textit{\ }

\smallskip

\item  If $\left\| f\right\| ^{\ast }_{\Lambda_2^p(uv)}:=\left( \int_{0}^{\infty
}\int_{0}^{\infty }(f_{yx}^{\ast \ast }(s,t))^pu(s)v(t)\,ds\,dt\right) ^{1/p})$, then  $\left\| \cdot \right\| ^{\ast }_{\Lambda_2^p(uv)}$ is a norm equivalent to $\left\| \cdot \right\| _{\Lambda_2^p(uv)}$.

\smallskip

\item  $S^{2}:L_{\rm{dec}}^{p}(\mathbb{R}^{2}_{+},uv)\rightarrow L^{p}(\mathbb{R}_{+}
^{2},uv)$ is bounded.

\smallskip

\item  $f\rightarrow f_{yx}^{\ast \ast }$ is bounded from $L_{\rm{dec}}^{p}(
\mathbb{R}^{2}_{+},uv)$ to $L^{p}(\mathbb{R}^{2}_{+},uv).$

\smallskip

\item  The norm $\left\| f\right\| ^{(2)}_{\Lambda_2^p(uv)}:=\left( \iint_{\mathbb{R} 
^{2}_+}(S^{2}f_{yx}^{\ast }(s,t))^{p}u(s)v(t)\,ds\,dt\right) ^{1/p}=\left\| f^{\ast \ast
}\right\|_{L^{p}(\mathbb{R}^{2}_{+},uv)}$ is equivalent to $\left\| f\right\|_{\Lambda_2^p(uv)} .$
\end{enumerate}
\end{theorem}

\noindent{\bf Proof:} 
(1.$\Longrightarrow $2.) Let $Tf(x,y)=f(x)g(y),$ where $g\in \Lambda ^{p}(v),$ $ 
\left\| g\right\| _{\Lambda ^{p}(v)}=1$ is a fixed function$.$ It is easy to
see that if $f\in \Lambda ^{p}(u)$, then $Tf\in \Lambda _{2}^{p}(uv)$ and $ 
\left\| Tf\right\|_{\Lambda_2^p(uv)} =\left\| f\right\| _{\Lambda ^{p}(u)},$ which shows that $\Lambda ^{p}(u)$ embeds isometrically into $\Lambda _{2}^{p}(uv)$. Since $\Lambda
_{2}^{p}(uv)$ is normable\textit{, }there exists a norm $\left\| \cdot
\right\| _{1}$ which is equivalent to $\left\| \cdot \right\|_{\Lambda_2^p(uv)} $ as defined
in (\ref{Lor2dim}). For a function $f\in \Lambda ^{p}(u),$ we define now $ 
\left\| f\right\| _{1}^{\ast }:=\left\| Tf\right\| _{1}.$ Since $\left\|
\cdot \right\| _{1}$ is a norm on $\Lambda _{2}^{p}(uv)$, it is clear that $ 
\left\|\cdot\right\| _{1}^{\ast }$ is a norm on $\Lambda ^{p}(u)$, which is equivalent
to $\left\| \cdot\right\| _{\Lambda ^{p}(u)}.$ Hence $u\in B_{p}.$ Similarly we
can prove that $v\in B_{p}.$

(2.$\Longleftrightarrow $3.)  This equivalence is a direct consequence of Proposition 
\ref{Lmix}.

(2.$\Longrightarrow $4.)  By Proposition \ref{properties}  we have
that 
\begin{equation}
\left\| f\right\|_{\Lambda _{2}^{p}(uv)}  =\left( \int_{0}^{\infty }\int_{0}^{\infty
}(f_{yx}^{\ast }(s,t))^{p}u(s)v(t)\,ds\,dt\right) ^{1/p} \leq \left\| f\right\|^{\ast }_{\Lambda _{2}^{p}(uv)}. \label{n1} 
\end{equation}

Since $v\in B_{p}$ and for any $s>0,\frac{1}{s}\int_{0}^{s}f_{yx}^{\ast
}(\sigma ,\cdot )\,d\sigma $ is a decreasing function, the boundedness of the
Hardy operator on $\Lambda^p(v)$, (see  \cite{AM}), gives 
\begin{equation}
\int_{0}^{\infty }\left( \frac{1}{t}\int_{0}^{t}\bigg(\frac{1}{s} 
\int_{0}^{s}f_{yx}^{\ast }(\sigma ,\tau)\,d\sigma \bigg) \,d\tau \right) 
^{p}v(t)\,dt\leq C_{1}^{p}\int_{0}^{\infty }\left( \frac{1}{s} 
\int_{0}^{s}f_{yx}^{\ast }(\sigma ,t)\,d\sigma \right) ^{p}v(t)\,dt,\text{ }
\label{H1}
\end{equation}
where $C_{1}>0,$ $s>0.$ Multiplying   (\ref{H1}) by $u(s)$, and
integrating with respect to $s$ gives, by the same arguments as above, that
there exists $C_{2}>0$ such that: 
\begin{align}
 \int_{0}^{\infty }\int_{0}^{\infty }(f_{yx}^{\ast \ast
}(s,t))^{p}u(s)v(t)\,ds\,dt  \notag 
 &=\int_{0}^{\infty }\int_{0}^{\infty }\left( \frac{1}{t}\int_{0}^{t}\left( 
\frac{1}{s}\int_{0}^{s}f_{yx}^{\ast }(\sigma ,\tau )\,d\sigma \right) \,d\tau  
\right) ^{p}u(s)v(t)\,ds\,dt  \notag \\
& \leq C_{1}^{p}\ \int_{0}^{\infty }\int_{0}^{\infty }\left( \frac{1}{s} 
\int_{0}^{s}f_{yx}^{\ast }(\sigma ,t)\,d\sigma \right)^{p}v(t)u(s)\,ds\,dt  \notag
\\
 &\leq (C_{1}^{\ }C_{2})^{p}\int_{0}^{\infty }\left( \int_{0}^{\infty
}(f_{yx}^{\ast }(s,t))^{p}v(t)\,dt\right) u(s)\,ds  \notag \\
 &=(C_{1}^{\ }C_{2})^{p}\left\| f\right\|_{\Lambda _{2}^{p}(uv)} ^{p}.  \label{n2}
\end{align}
Combining now (\ref{n1}) and (\ref{n2}) we get 
\begin{equation*}
\left\| f\right\|_{\Lambda _{2}^{p}(uv)}   \leq \left\| f\right\| _{\Lambda _{2}^{p}(uv)}  ^{\ast }\leq C\left\| f\right\|_{\Lambda _{2}^{p}(uv)}  .
\end{equation*}

On the other hand, by Proposition \ref{properties},  $\left\| f\right\|
^{\ast }_{\Lambda _{2}^{p}(uv)} $ is \ a norm.

(4.$\Longrightarrow $1.) The fact that $\left\| \cdot \right\|_{\Lambda _{2}^{p}(uv)} ^{\ast }$ is an
equivalent norm on $\Lambda _{2}^{p}(uv)$ implies that $\Lambda _{2}^{p}(uv)$
is normable.

(4.$\Longleftrightarrow $5.$\Longleftrightarrow $6.)  This is obvious according to
Lemma \ref{LH} and the definitions of $S^{2}$ and $S_{2,1}.$

(7.$\Longrightarrow $5.)  This implication is trivial.

(5.$\Longrightarrow $7.) By hypothesis we have that $\left\| f\right\|_{\Lambda _{2}^{p}(uv)}^{(2)}\le C
\left\| f\right\|_{\Lambda _{2}^{p}(uv)} $, for some $C>0$. On the other hand by Proposition \ref{properties} we
have $\left\| f\right\|_{\Lambda _{2}^{p}(uv)} \leq \left\| f\right\| _{\Lambda _{2}^{p}(uv)}^{\ast }\leq \left\|
f\right\|_{\Lambda _{2}^{p}(uv)} ^{(2)}$ and this completes the proof.
$\hfill\Box$\bigbreak

Now we will compare the spaces $\Lambda_2^p(uv)$ and $\Lambda^p(u)[\Lambda^p(v)]$.
The next theorem gives an embedding result between the two-dimensional
Lorentz space $\Lambda _{2}^{p}(uv )\ $ and the mixed norm space $ 
\Lambda  ^{p}(u)[\Lambda^{p}(v)].$

\begin{theorem}\label{incdec}
Let   $u,v:\mathbb{R}
_{+}\mathbb{\rightarrow R}_{+}$ be two weights, and assume that $u$ is a decreasing function. Then \ $\Lambda
_{2}^{p}(uv)\subset \Lambda  ^{p}(u)[\Lambda^{p}(v)].$
\end{theorem}

\noindent{\bf Proof:}
As in the proof of Proposition \ref{properties}, and using (\ref{f1}), we have 
\begin{align*}
 \int_{0}^{\infty }\left( \int_{0}^{\infty }(f_{y}^{\ast }(\cdot
,t))^pv(t)\,dt\right) _{x}^{\ast }(s)u(s)\,ds 
& =\sup_{\rho }\int_{\mathbb R }\left( \int_{0}^{\infty }f_{y}^{\ast
p}(x,t)v (t)\,dt\right) u(\rho (x))\,dx \\
& =\sup_{\rho }\int_{0}^{\infty }\left( \int_{\mathbb R }f_{y}^{\ast
p}(x,t) u(\rho (x))\,dx\right) v (t)\,dt \\
& \leq \int_{0}^{\infty }\left( \sup_{\rho }\int_{\mathbb R }f_{y}^{\ast
p}(x,t) u(\rho (x))\,dx\right) v (t)\,dt \\
& =\int_{0}^{\infty }\int_{0}^{\infty }f_{yx}^{\ast
p}(s,t)u(s)v (t)\,ds\,dt.\qquad\qquad\qquad\Box
\end{align*}
\bigbreak

\begin{remark}\label{noincl}{\rm
\label{Rem}  (i) Although $\Lambda _{2}^{p}(uv)\ $ and $\Lambda^{p}(u)[\Lambda^{p}(v)]$ are both Banach spaces if and only if $u$ and $v$ satisfy the $ 
B_{p}$ condition, the two spaces are not  equal. We will prove this fact by
means of the following example. Let $u(s)=\chi _{\lbrack 0,1]}(s)$, $v\equiv
1\ $and 
\begin{equation*}
f(x,y)= \begin{cases} {a_{[x]}\chi _{([x],[x]+1)\times (0,[x]+1);}  } & \text{ $x,y\geq 0$}, \\
      0& \text{otherwise},
\end{cases}
\end{equation*}
where $a_k=1/(1+k)^{1/p}$, and  $[x]$  
is the integer part of $x$.

We have 
\begin{equation*}
f_{y}^{\ast }(x,t)=a_{[x]}\chi _{(0,1+[x])}(t)
\end{equation*}
and 
\begin{equation*}
f_{yx}^{\ast }(s,t)=(f_{y}^{\ast }(\cdot ,t))_{x}^{\ast
}(s)=\sum_{k=0}^{\infty }a_{[t]+k}\chi _{(k,k+1)}(s).
\end{equation*}
Hence, 
\begin{align*}
\left\| f\right\| _{\Lambda^{p}(u)[\Lambda^{p}(v)]}& =\ \left( \int_{0}^{1}\left( \int_{0}^{\infty
}f_{y}^{\ast p}(\cdot ,t)\,dt\right) _{x}^{\ast }(s)\,ds\right) ^{1/p}\  \\
& =\ \ \left( \int_{0}^{1}\left( \int_{0}^{1+[x]}a_{[x]}^{p}\,dt\right)
_{x}^{\ast }(s)\,ds\right) ^{1/p}=\left( \int_{0}^{1}\left(
a_{[x]}^{p}(1+[x])\right) _{x}^{\ast }(s)\,ds\right) ^{1/p} \\
& = \left( \int_{0}^{1}1\,ds\right) ^{1/p}=1,
\end{align*}
and  
\begin{align*}
\left\| f\right\|_{\Lambda _{2}^{p}(uv)}  & =\left( \int_{0}^{\infty }\int_{0}^{1}f_{yx}^{\ast
p}(s,t)\,ds\,dt\right) ^{1/p} \\
& =\left( \int_{0}^{\infty }\int_{0}^{1}a_{[t]}^{p}\,ds\,dt\right) ^{1/p}\
=\left( \int_{0}^{\infty }a_{[t]}^{p}\,dt\right) ^{1/p} \\
& =\left( \sum_{k=0}^{\infty }a_{k}^{p}\right) ^{1/p}=\left( \sum_{k=0}^{\infty }\frac1{1+k}\right) ^{1/p}=\infty.
\end{align*}
Hence,   we see that $f\in \Lambda^{p}(u)[\Lambda^{p}(v)]\setminus \Lambda _{2}^{p}(uv)$.

\bigbreak

(ii) 
Using Theorem \ref{incdec} and (i) we have that, if $u=\chi_{[0,1]}$ and $v=1$, then $\Lambda _{2}^{p}(uv)$ is strictly contained in $\Lambda^{p}(u)[\Lambda^{p}(v)]$. We will now prove that, with the same weights,  $\Lambda _{2}^{p}(uv)\not\subset\Lambda^{p}(v)[\Lambda^{p}(u)]$: In fact, take $f(x,y)=\sum_{k=0}^{\infty} a_k\chi_{([k,k+1]\times[k,k+1])}(x,y)$, where $\{a_k\}_k$ is a decreasing sequence. Then, $f^*_{yx}(s,t)=\sum_{k=0}^{\infty}a_k\chi_{([k,k+1]\times[0,1])}(s,t)$, and $\Vert f\Vert_{\Lambda _{2}^{p}(uv)}=a_0$. However, $f^*_x(s,y)=a_{[y]}\chi_{[0,1]}(s,y)$, and hence, $\Vert f\Vert_{\Lambda^{p}(v)[\Lambda^{p}(u)]}=\big(\sum_{k=0}^{\infty}a_k^p\big)^{1/p}=\infty$, if we take, for example, $a_k=1/(k+1)^{1/p}$. Therefore, we have shown that neither of the three spaces $\Lambda _{2}^{p}(uv)$, $\Lambda^{p}(u)[\Lambda^{p}(v)]$, and $\Lambda^{p}(v)[\Lambda^{p}(u)]$ are equal.
}
\end{remark}

\bigbreak

\begin{remark}\label{cotrau}{\rm
Theorem \ref{mainth} is false for $p=1$ (as in the one-dimensional case): In fact, take $u=1$, $v=1$. Then $\Lambda _{2}^{1}(uv)=L^1(\mathbb R^n)$, which is Banach, but $1\notin B_1$. 
\bigbreak

\noindent
We could also try to consider the case of general weights in Theorem \ref{mainth}. For this, we define the following classes (see \cite{BS} for the definition of the weak-type spaces):
}\end{remark}

\begin{definition}
For $p\ge 1$ we say that $w\in B^{(2)}_p$, if $S^{2}:L_{\rm{dec}}^{p}(\mathbb{R}^{2}_{+},w)\rightarrow L^{p}(\mathbb{R}_{+}
^{2},w)$ is bounded, and $w\in B^{(2)}_{p,\infty}$, if $S^{2}:L_{\rm{dec}}^{p}(\mathbb{R}^{2}_{+},w)\rightarrow L^{p,\infty}(\mathbb{R}_{+}
^{2},w)$ is bounded.
\end{definition}

\begin{definition}
We say that a set $D\subset \mathbb{R}_{+}^{2}$ is decreasing (and write $ 
D\in \Delta _{d}$) if the function $\chi _{D}$ is decreasing in each
variable.
\end{definition}

\noindent
We will now show several properties of these weights, and prove some extensions of Theorem~\ref{mainth}.

\begin{theorem}\label{prodw} Assume that $p\ge 1$.

\begin{enumerate}

\item If $w\in B^{(2)}_p$, then $\Lambda _{2}^{p}(w)$ is normable.

\item $w\in B^{(2)}_1$ if and only if $\sup_{D\in\Delta_d }\displaystyle\frac{\int_{\mathbb R^2_+}S^{(2)}(\chi_D)(s,t)w(s,t)\,ds\,dt}{\int_Dw(s,t)\,ds\,dt}<\infty$.

\item If $w(s,t)=u(s)v(t)$, then 
$$
\sup_{D\in\Delta_d  }\displaystyle\frac{\int_{\mathbb R^2_+}S^{(2)}(\chi_D)(s,t)w(s,t)\,ds\,dt}{\int_Dw(s,t)\,ds\,dt}=\left(1+\sup_{a>0}\frac{a\int_a^{\infty}\frac{u(s)}s\,ds}{\int_0^au(s)\,ds}\right)\left(1+\sup_{b>0}\frac{b\int_b^{\infty}\frac{v(t)}t\,dt}{\int_0^bv(t)\,dt}\right).
$$

\item If $w(s,t)=u(s)v(t)$, then $w\in B^{(2)}_p$ if and only if $u,v\in B_p$.

\item $1\notin B^{(2)}_{1,\infty}$, although $1\in B_{1,\infty}$.
\end{enumerate}
\end{theorem}

\noindent{\bf Proof:} 
(1.) Similarly as in Theorem \ref{mainth}, we have that 
$$
\left\| f\right\| _{\Lambda_2^p(w)}\le\left\| f\right\| ^{*}_{\Lambda_2^p(w)}=\bigg(\int_0^{\infty}\int_0^{\infty}(f^{**}_{yx}(s,t))^pu(s)v(t)\,dsdt\bigg)^{1/p}\le\left\| f\right\| ^{(2)}_{\Lambda_2^p(w)},
$$
and hence, since $w\in B^{(2)}_p$, 
$$
\left\| f\right\| ^{(2)}_{\Lambda_2^p(w)}\le C\left\| f\right\| _{\Lambda_2^p(w)},
$$
showing that $\left\| f\right\| _{\Lambda_2^p(w)}$ is equivalent to $\left\| f\right\| ^{*}_{\Lambda_2^p(w)}$. By Proposition \ref{properties} (2), it is clear that $||\cdot||_{\Lambda_2^p(w)}^*$ is a norm.

\medskip 

(2.) This is a special case of \cite[Theorem 2.2 (c)]{BaPeSo1}.

\medskip

(3.) If we write $\tilde u(\sigma)= \int_{\sigma}^{\infty}\frac{u(s)}s\,ds$, and similarly $\tilde v$, then
$$
 \int_{\mathbb R^2_+}S^{(2)}(\chi_D)(s,t)u(s)v(t)\,ds\,dt = \int_D \tilde u(\sigma)\tilde v(\tau)\,d\sigma\,d\tau,
 $$
and hence, using  \cite[Theorem 2.5]{BaPeSo1}
\begin{eqnarray*}
&&\sup_{D\in\Delta_d  }\displaystyle\frac{\int_{\mathbb R^2_+}S^{(2)}(\chi_D)(s,t)u(s)v(t)\,ds\,dt}{\int_Du(s)v(t)\,ds\,dt}\\&=&\sup_{D \in\Delta_d }\displaystyle\frac{\int_D \tilde u(\sigma)\tilde v(\tau)\,d\sigma\,d\tau}{\int_Du(s)v(t)\,ds\,dt}=
\sup_{a,b>0}\frac{\big(\int_0^a\tilde u(\sigma)\,d\sigma\big)\big(\int_0^b\tilde v(\tau)\,d\tau\big)}{\big(\int_0^a  u(s)\,ds\big)\big(\int_0^b v(t)\,dt\big)}\\
&=&
\bigg(\sup_{a>0}\frac{\int_0^a  u(s)\,ds+a\int_a^{\infty}\frac{u(s)}s\,ds}{\int_0^a  u(s)\,ds}
\bigg)\bigg(\sup_{b>0}\frac{\int_0^b  v(t)\,dt+b\int_b^{\infty}\frac{v(t)}t\,dt}{\int_0^b  v(t)\,dt}\bigg).
\end{eqnarray*}

\medskip

(4.) If $p=1$, then the conclusion follows immediately by using (2.) and (3.). If $p>1$,   then by using (1.) we have that if  $w(s,t)=u(s)v(t)\in B_p^{(2)}$, then $\Lambda _{2}^{p}(w)$ is a Banach space, and by Theorem \ref{mainth} we conclude that $u,v\in B_p$. Conversely, if $u,v\in B_p$, then it suffices to observe that if $f(\sigma,\tau)$ is a decreasing function, then
$$
S^{(2)}f(s,t)=S(S(f_{\tau}(s))(t),
$$
which is a composition of the one-dimensional Hardy operators acting on decreasing functions, and hence, for product weights, if $S$ is bounded on both $L^p_{\rm dec}(u)$ and $L^p_{\rm dec}(v)$, then $S^{(2)}$ is bounded on $L^p_{\rm dec}(uv)$.

\medskip

(5.) It is well known that $1\in B_{1,\infty}$. To show that $1\notin B^{(2)}_{1,\infty}$, we consider the function $f=\chi_{[0,1]\times[0,1]}$. Then, 
$$
S^{(2)}f(s,t)=\frac1{st}\min\{s,1\}\min\{t,1\},
$$
and hence, for $0<\lambda<1$, 
$$
\lambda|\{(s,t):\, S^{(2)}f(s,t)>\lambda\}|\ge\lambda\int_1^{1/\lambda}\frac1{\lambda s}\,ds-\lambda (\frac1\lambda-1)=\log\frac1{\lambda}+\lambda-1.
$$
Therefore, $S^{(2)}f\notin L^{1,\infty}$.
$\hfill\Box$\bigbreak

\begin{remark}{\rm (i) 
 Contrary to what one could expect looking at Theorem \ref{mainth} and the one-dimensional case (\cite{CaSoA}), Theorem \ref{prodw} shows that $\Lambda^1_2(w)$ can be a Banach space, but $w$ need not be in $B^{(2)}_{1,\infty}$, even in the product case (just take $w=1$; see also Remark \ref{cotrau}).

\medskip

(ii) There are other related properties, like linearity or quasinormability, which are also of interest in this setting. For the Lorentz spaces $\Lambda^p(u)$, quasinormability was characterized in  \cite{H} (see also  \cite{CS} and \cite{KM}), and for the multidimensional version, $\Lambda^p_2(w)$, in \cite{BaPeSo}. However, linearity (when is it a vector space?) is only known for  $\Lambda^p(u)$ (see \cite{CKMP}). }
\end{remark}

\bigbreak

\section{\protect\bigskip Embedding results}\label{emth}

This section deals with embedding results between the different kind of Lorentz spaces introduced before. We will start by\ proving a theorem  which
characterizes the embeddings between the classical Lorentz spaces and those
defined in (\ref{Lor2dim}), for the case $p\le q$.

\begin{theorem}
\label{incl}Let $0<p\leq q<\infty ,u:\mathbb{R}_{+}\mathbb{\rightarrow R} 
_{+} $ and $w:\mathbb{R}_{+}^{2}\mathbb{\rightarrow R}_{+}$ be two weight
functions. The following conditions are equivalent:

\begin{enumerate}
\item  $\Lambda ^{p}(\mathbb{R}^2,u)\subset \Lambda _{2}^{q}(w)$.
\item  $C=\sup_{D\in \Delta _{d}}\displaystyle\frac{\left( \int_{D}w(s,t)\,ds\,dt\right) ^{1/q}}{ 
\left( \int_{0}^{\left| D\right| }u(\theta)\,d\theta\right) ^{1/p}}<\infty $, and $C$ is the
best constant for the embedding.
\end{enumerate}
\end{theorem}

\noindent{\bf Proof:} 
By evaluating the inequality of the embedding in (1.), for the decreasing function $f=\chi
_{D}$, we get (2.).

To prove that (2.) implies (1.), by the easy observation that $\ f_{yx}^{\ast }$
is equimeasurable with $f$, we may assume that $f$ is a decreasing function in each variable. We will follow the ideas and  techniques   from 
\cite{BaPeSo1}. In particular,  we will use the embedding $L^{1}_{\rm{dec}}(t^{p-1})\subset
L^{q/p}(t^{q-1}),$ and the remark   that $D_{\lambda}=\{(s,t):f(s,t)>\lambda\}$ is a decreasing
set:
\begin{align*}
\left( \iint_{\mathbb{R}^{2}_+}\ f^{q}(s,t)w(s,t)\,ds\,dt\right) ^{1/q}&
=q^{1/q}\left( \int_{0}^{\infty }\lambda^{q-1}\bigg(\int_{D_{\lambda}}w(s,t)\,ds\,dt\right)\,d\lambda\bigg) ^{1/q}
\\
& \leq Cq^{1/q}\left( \int_{0}^{\infty }\lambda^{q-1}\left( \int_{0}^{\left|
D_{\lambda}\right| }u(\theta )\,d\theta \right) ^{q/p}\,d\lambda\right) ^{1/q} \\
& \leq Cq^{1/q}\left( \frac{q}{p}\right) ^{-1/q}p^{1/p-1/q}\left(
\int_{0}^{\infty }\lambda^{p-1}\bigg(\int_{0}^{\left| D_{\lambda}\right| }u(\theta )\,d\theta\bigg)\,d\lambda
\right) ^{1/p} \\
& =Cp^{1/p}\left( \int_{0}^{\infty }\lambda^{p-1}\bigg(\int_{0}^{\left| D_{\lambda}\right| }u(\theta )\,d\theta\bigg)\,d\lambda \right) ^{1/p} \\
& =C\ \left( \int_{0}^{\infty }f^{\ast p}(\theta )u(\theta )\,d\theta \right)
^{1/p}.
\end{align*}
Hence we also obtain that $C$ is the best constant. 
$\hfill\Box$\bigbreak

\begin{theorem}
\label{incl1}Let $0<p\leq q<\infty ,$ $u:\mathbb{R}_{+}\mathbb{\rightarrow R} 
_{+}$ and $w:\mathbb{R}_{+}^{2}\mathbb{\rightarrow R}_{+}$ be two weight
functions. The following conditions are equivalent:

\begin{enumerate}
\item  $\Lambda _{2}^{p}(w)\subset \Lambda ^{q}(\mathbb{R}^2,u)$.

\smallskip

\item  $C=\sup_{D\in \Delta _{d}}\displaystyle\frac{\left( \int_{0}^{\left| D\right|
}u(\theta)\,d\theta\right) ^{1/q}}{\left( \int_{D}w(s,t)\,ds\,dt\right) ^{1/p}}<\infty $, and $C$ is the
best constant for the embedding.\end{enumerate}
\end{theorem}

\noindent{\bf Proof:} 
By evaluating the inequality of the embedding in (1.), for the decreasing function $f=\chi
_{D}$, we get (2.)

To prove that (2.) implies (1.), it is enough to consider the inequality only
for decreasing functions in each variable and to apply the same technique as
in the proof of Theorem \ref{incl}. 
$\hfill\Box$\bigbreak

\begin{remark}{\rm
If $p=q$ we get from   Theorems \ref{incl} and  \ref{incl1} that $\Lambda
_{2}^{p}(w)=\Lambda ^{p}(\mathbb{R}^2,u)$ if and only if $\int_{0}^{\left| D\right|
}u(\theta)\,d\theta\approx \int_{D}w(s,t)\,ds\,dt,$ for all decreasing sets $D$.  Theorems \ref{incl} and \ref{incl1} generalize 
\cite[Theorem 1]{Y}.}
\end{remark}

We now consider the remaining case $p>q$. The proofs of these results follow the same arguments used in \cite{BPSt} and \cite{JL}, with small modifications, and hence we will omit them. We will use the following notations: A covering family of $\mathbb R_+^2$ is an increasing family of decreasing sets $\{D_k\}_k$, such that $\cup_kD_k=\mathbb R_+^2$. In this case, we set $\Delta_k=D_{k+1}\setminus D_k$. If $f$ is a function in $\mathbb R_+^2$, decreasing in each variable, $D_{f,t}=\{f>t\}$. Recall also that if $w$ is a weight in $\mathbb R_+^2$ and $E$ is a measurable set, then $w(E)=\int_Ew(s,t)\,ds\,dt$, and if $u$ is a weight in $\mathbb{R}_{+}$, we write $U(t)=\int_0^tu(s)\,ds$.

\begin{theorem}
\label{jl1} Let $0<q<p<\infty ,$ $u:\mathbb{R}_{+}\mathbb{\rightarrow R} 
_{+}$ and $w:\mathbb{R}_{+}^{2}\mathbb{\rightarrow R}_{+}$ be two weight
functions. Let $r=pq/(p-q)$. The following conditions  are equivalent:

\begin{enumerate}
\item  $\Lambda ^{p}(\mathbb{R}^2,u)\subset \Lambda_{2} ^{q}(w)$

\smallskip

\item  There exists a constant $C>0$ such that for all covering families $\{D_k\}_k$:
$$
\int_0^1\bigg(\sum_k\bigg(\frac{w(D_k)+w(\Delta_k)t}{U(|D_k|)+(U(|D_{k+1}|)-U(|D_k|))t}\bigg)^{r/p}w(\Delta_k)\bigg)\,dt\le C.
$$

\item  There exists a constant $C>0$ such that for all covering families $\{D_k\}_k$:
$$
\sum_k\big(w(\Delta_k)^{r/q}\,U(|D_{k+1}|)^{-r/p}\big)\le C.
$$

\item  There exists a constant $C>0$ such that for all two-dimensional functions   $f$ in $\mathbb R_+^2$, decreasing in each variable:
$$
\int_0^{\infty}U(|D_{f,t}|)^{-r/p}\,d\big(-w(D_{f,t})^{r/q}\big)\le C.
$$

\end{enumerate}
\end{theorem}
We also state  the converse embedding:

\begin{theorem}
\label{jl2} Let $0<q<p<\infty ,$ $u:\mathbb{R}_{+}\mathbb{\rightarrow R} 
_{+}$ and $w:\mathbb{R}_{+}^{2}\mathbb{\rightarrow R}_{+}$ be two weight
functions. Let $r=pq/(p-q)$. The following conditions are equivalent:

\begin{enumerate}
\item  $\Lambda_{2} ^{p}(w)\subset \Lambda ^{q}(\mathbb{R}^2,u)$.

\smallskip

\item  There exists a constant $C>0$ such that for all covering families $\{D_k\}_k$:
$$
\int_0^1\bigg(\sum_k\bigg(\frac{U(|D_k|)+(U(|D_{k+1}|)-U(|D_k|))t}{w(D_k)+w(\Delta_k)t}\bigg)^{r/p}(U(|D_{k+1}|)-U(|D_k|))\bigg)\,dt\le C.
$$

\item  There exists a constant $C>0$ such that for all covering families $\{D_k\}_k$:
$$
\sum_k\big(w(D_{k+1})^{-r/p}\,(U(|D_{k+1}|)-U(|D_k|))^{r/q}\big)\le C.
$$

\item  There exists a constant $C>0$ such that for all two-dimensional functions   $f$ in $\mathbb R_+^2$, decreasing in each variable:
$$
\int_0^{\infty}w(D_{f,t})^{-r/p}\,d\big(-U(|D_{f,t}|)^{r/q}\big)\le C.
$$

\end{enumerate}
\end{theorem}

\bigbreak

\address{
\noindent
Sorina Barza\\ Dept. of Mathematics, Physics and Eng. Sciences
\\ Karlstad University\\ SE-65188 Karlstad, SWEDEN\ \ \ {\sl E-mail:} 
  {\tt sorina.barza@kau.se}

\bigbreak
\noindent
Anna Kami\'nska\\ Dept. of Mathematical Sciences\\ The University of Memphis\\ Memphis, TN,
 38152-3240, USA\ \ \ {\sl E-mail:} {\tt  kaminska@memphis.edu}

\medskip
\noindent
Lars-Erik Persson\\ Dept. of Mathematics
\\ Lule\aa\ University of Technology\\ SE-97187 Lule\aa,
 SWEDEN\ \ \ {\sl E-mail:} 
 {\tt larserik@sm.luth.se}

\medskip
\noindent
Javier Soria\\ Dept. Appl. Math. and Analysis
\\ University of Barcelona\\ E-08071 Barcelona,
 SPAIN\ \ \ {\sl E-mail:} 
 {\tt soria@mat.ub.es}
}

\end{document}